\documentclass[11pt,intlimits]{amsart}
\usepackage{amssymb,amsmath,amsthm,upref,cite}

\numberwithin{equation}{section}

\theoremstyle{plain}

\theoremstyle{remark}

\usepackage{graphicx}

\usepackage{amsmath,amssymb,eucal,graphicx}
\newcommand{\be}{\begin{equation}}
\newcommand{\ee}{\end{equation}}
\newcommand{\la}{\label}

\newcommand{\beq}{\begin{equation}}
\newcommand{\eeq}{\end{equation}}
\newcommand{\dd}{{d}}
\newcommand{\ii}{{i}}

\begin{document}

\title[Limits of measure redistribution processes]{Universal limits of nonlinear measure redistribution processes and their applications}

\author[Razvan Teodorescu]{Razvan Teodorescu}
\address{Department of Mathematics \& Statistics, University of South Florida, Tampa, FL 33620-5700, USA}
\email{razvan@usf.edu}

\date{}


\maketitle

\begin{abstract}

Deriving the time evolution of a distribution of probability (or a probability density matrix) is a 
problem encountered frequently in a variety of situations: for physical time, it could be a kinetic reaction 
study, while identifying time with the number of computational steps gives a typical picture of algorithms 
routinely used in quantum impurity solvers, density functional theory, etc. Using a truncation scheme 
for the expansion of the exact quantity is necessary due to constraints of the numerical implementation. 
However, this leads in turn to serious complications such as the Fermion Sign Problem (essentially, 
density or weights will become negative). By integrating angular degrees of freedom and reducing the 
dynamics to the radial component,  the time evolution is reformulated as a nonlinear integral transform 
of the distribution function. A canonical decomposition into orthogonal polynomials leads back to the 
original sign problem, but using a characteristic-function representation allows to extract the asymptotic 
behavior, and gives an exact large-time limit, for many initial conditions, with guaranteed positivity. 

\end{abstract}


\section{Introduction}

The notion of coarse graining in statistical physical models (or field theory), 
introduced by  Migdal and Kadanoff \cite{Migdal1975}, \cite{Kadanoff1976}
is essential to many fundamental concepts and results, like the continuum 
limit of lattice models, or the universality of scaling behavior near a phase 
transition, to name two of the most celebrated consequences. 

From the perspective of mathematical statistics, the analysis of such coarse 
graining processes is straightforward due to the fact that the elementary operation of the process is $averaging$: at step $n$, we create a new random variable $X^{(n)}$ from two variables defined at step $n-1$ by $X^{(n)}= [X^{(n-1)}_1+X^{(n-1)}_2]/2$. Then by repeated application of this operation, the result after sufficiently many steps is simply given by the central limit theorem (CLT). 

In this work, we consider another class of coarse-graining processes, where at each step we 
take the absolute difference between variables, rather than their average. This is justified by a number of relevant physical problems, but also by classical issues from decision theory or economics. 
As in the case of standard coarse graining, there is a limiting distribution (in fact, a whole class) which will be reached after arbitrarily many steps. Unlike in the standard case, the particular limit is chosen from this class based on the asymptotic properties of the initial distribution (more precisely, the radius of convergence of its moment-generating function). This is an example of the {\emph {extreme  selection  criterion }} which characterizes other important stochastic processes, such as the Fisher-Kolmogorov evolution. 

The structure of this paper is the following: in the first section we present the difference coarse-graining procedure, as well as some of its realizations, and derive its universal long-time asymptotic behavior. 
In the second section, we consider some particular types of distributions which are useful examples for the general results. The last section is a discussion on possible applications of this new universal limiting behavior. 

\section{Difference coarse-graining processes and asymptotic limits}

\subsection{Excess redistribution and partial annihilation models}

In the excess process, one starts with random variables
$-\infty<y_i<\infty$ for $i=1,\ldots,N$. These variable are initially
distributed according to some prescribed distribution $P_0(y)$. Then,
two variables $y_1$ and $y_2$ are picked at random. If both are
positive or both are negative, nothing is done. If one is positive and
one is negative, then this variables are updated as follows
\begin{equation}
(y_1,y_2)\to (y_1+y_2,0).
\end{equation}
Thus, this process conserves the total sum $\sum_i y_i$. A 
related process was considered in \cite{Krapivsky-Ben-Naim}

Let $P(x,t)$ be the time-dependent probability density of the process. 
It evolves according to the rate equation
\begin{equation}
\label{P-eq-0}
\frac{\partial P(x)}{\partial t}=-\,c(t)P(x)+\int_{-\infty}^0 dy\, P(y)P(x-y).
\end{equation}
This equation holds for $x>0$ and a similar equation holds for
$x<0$. The integral has a convolution form. Here $c(t)$ is the total
density of non-zero variables.

One motivation for this process comes from economics. A positive $y$
represents wealth and a negative $y$ represents debt. The economy
evolves under conservative exchange of excess wealth. Another
motivation is the electric power-grid where power producing plants may
exchange excess or deficit capacity in response to power demands.

Let us assume that the initial distributions are symmetric
$P_0(y)=P_0(-y)$. Clearly, this property persists with
time. Therefore, the above process may be equivalently formulated by
considering the evolution of the variable $x=|y|$. Then, the excess
process becomes
\begin{equation}
\label{process-a}
(x_1,x_2)\to |x_1-x_2|.
\end{equation}

Let $c(t)$ be the remaining fraction of non-zero variables at time
$t$. This quantity satisfies
\begin{equation}
\label{c-eq}
\frac{dc}{dt}=-c^2.
\end{equation}
Solving this equation subject to the initial condition $c(0)=1$, we find
\begin{equation}
\label{c-sol}
c(t)=\frac{1}{1+t}.
\end{equation}

Next, we consider $P(x,t)$, the probability density. This quantity
evolves according to the integro-differential equation
\begin{equation}
\label{P-eq}
\frac{\partial P(x)}{\partial t}=-2\,c(t)P(x)+2\int_0^\infty dy\, P(y)P(x+y).
\end{equation}
The negative term accounts for loss of two interacting variables and
the gain term accounts for gain of one variable. Of course, the total
density is $c(t)=\int dx P(x)$. Integrating this equation, we recover
(\ref{c-eq}). Despite the simple nature of this equation, it is
challenging. For example, the first moment $M_1(t)=\int_0^\infty xP(x) dx
$ does not obey a closed equation.

Consider the normalized distribution
\begin{equation} \label{normalization}
p(x,t)=c^{-1}P(x,t),
\end{equation}
so that $\int dx\, p(x,t)=1$. It satisfies the evolution equation
\begin{equation}
\label{p-eq}
\frac{\partial p(x)}{\partial t}=-c\,p(x)+2c\int_0^\infty dy\, p(y)\,p(x+y).
\end{equation}
Now, we introduce the new time variable
\begin{equation}
\label{tau-def}
\tau=\ln (1+t).
\end{equation}
The master equation becomes
\begin{equation}
\label{pp-eq}
\frac{\partial p(x)}{\partial \tau}=-p(x)+2\int_0^\infty dy\, p(y)p(x+y).
\end{equation}
This master equation represents the dual process $(x_1,x_2)\to
(|x_1-x_2|,|x_1-x_2|)$.  This process occurs with rate $1/2$. In this
formulation, the number of particles is conserved $\int dx\,
p(x)=1$. 

In the following section we explore the general properties of the 
long-time behavior of the probability density function. We prove that there
exists a whole class of exact solutions, and quantify the convergence 
rate of an arbitrary distribution towards the appropriate infinite-time 
limit. 

\section{Convergence properties of the difference coarse graining}

\subsection{Steady-state solutions}

At the steady state, the probability density satisfies
\begin{equation} \la{steady0}
p(x)=2\int_0^\infty\, dy\, p(y)\,p(x+y).
\end{equation}
Performing the Taylor expansion of $p(x+y)$ with respect to variable $y$ and integrating, 
we obtain 
\be \la{steady1}
p(x) = 2 \sum_{n=0}^\infty \frac{M_n}{n!} \partial_x^n p(x) = 2 m(\partial_x) p(x),
\ee
where $m(t) = \sum_{n=0}^\infty \frac{M_n t^n}{n!}$ is the moment generating function
of $p(x)$, and $M_n$ is its non-centered moment of order $n$. Equation (\ref{steady1})
can also be written in the eigenvalue form
\be
m^{-1}(\partial_x) p(x) = 2 p(x), 
\ee
with $m^{-1}$ understood as the corresponding inverse pseudo-differential operator. 
In this form, it is very simple to notice that (\ref{steady1}) has 
a continuos family of solutions, exponential distributions $p(x) = \lambda 
e^{-\lambda x}$, $\lambda > 0$, for which the inverse moment generating function
is $m^{-1}_\lambda(t) = 1 -  t \lambda^{-1}$:
\be \la{solutions}
[1-\lambda^{-1} \partial_x] e^{-\lambda x} = 2 e^{-\lambda x}. 
\ee
Identity (\ref{solutions}) shows that there exists a class of solutions indexed by 
the continuous parameter $\lambda > 0$
\begin{equation} 
\label{steady-sol}
p(x)=\lambda\,e^{-\lambda\, x}.
\end{equation}
The fact that this coarse-graining process evolves towards a universal limit (in 
functional sense) bears some resemblance to the Central Limit Theorem. In that 
case, the limit (properly rescaled) is determined by the second  moment of the 
initial distribution. We will show that in the case of difference processes, it is not 
a single centered moment which determines the limit, but rather the convergence 
properties of the whole momen-generating function. In the next section, we investigate 
how initial conditions determine the particular steady-state of type (\ref{steady-sol}) 
realized by a particular process, as well as  the convergence properties of the evolution 
towards the solution.

\subsection{Asymptotic analysis}

In order to establish the steady-state reached by a difference coarse-graining 
process, as well as to characterize its convergence, we reformulate the evolution
(\ref{pp-eq}) through characteristic functions, since convergence in distribution 
(weak convergence) is equivalent to point-wise convergence of characteristic
functions,  $\phi(z, \tau)=\mathbb{E}[e^{izx}]$.

Define the charateristic function of distribution $P(x,\tau)$ as
\beq \label{definition}
\phi(z,\tau) = \int_0^\infty e^{\ii zx} P(x, \tau) \dd x,
\eeq
and obtain for the evolution equation
\beq
P(x, \tau) + \dot{P}(x,\tau) = 2 \int_0^\infty P(y, \tau)P(x+y, \tau) \dd y
\eeq
the form
$$
\phi(z, \tau) + \dot{\phi}(z, \tau) = 
2\int_0^\infty \int_0^\infty \int_0^\infty P(y, \tau)P(w, \tau) e^{\ii zx} \delta(w - x-y) \dd y \dd w \dd x.
$$
Use the standard representation for Dirac's distribution and obtain
$$
\phi(z, \tau) + \dot{\phi}(z, \tau) = 
2\int_0^\infty \int_0^\infty \int_0^\infty \int_{-\infty}^\infty P(y, \tau)P(w, \tau) 
e^{\ii [zx + k(w-y-x)]} \dd y \dd w \dd x \frac{\dd k}{2\pi}.
$$
Perform the integration over $x$ using the limit
\beq
\int_0^\infty e^{\ii x(z-k)} \dd x = \lim_{\epsilon \to 0^+} \int_0^\infty e^{\ii x(z-k + \ii \epsilon)} \dd x 
= \lim_{\epsilon \to 0^+}\frac{1}{\ii(z-k+\ii \epsilon)}.
\eeq
Now integrate over $y, w$ and get
\beq
\phi(z, \tau) + \dot{\phi}(z, \tau) = 
\lim_{\epsilon \to 0^+} \frac{1}{2 \pi \ii}\int_{-\infty}^\infty \frac{2\phi(-k, \tau)}{z-k +\ii \epsilon} \phi(k, \tau) \dd k, 
\eeq
which can be written as the Cauchy integral over the boundary of the lower half-plane $\Lambda$, in clockwise direction:
\beq
\phi(z, \tau) + \dot{\phi}(z, \tau) = 
\lim_{\epsilon \to 0^+} \oint_{\Lambda} \frac{2\phi(-k, \tau)}{z-k +\ii \epsilon} \phi(k, \tau) \frac{\dd k}{2 \pi \ii}. 
\eeq

Equation (\ref{pp-eq}) becomes:
\begin{equation} \label{character}
\dot{\phi}(s, \tau) = - \phi(s, \tau) + \frac{1}{2\pi i}\oint_{\Lambda}
\frac{\phi(-z, \tau)}{s-z+i \epsilon}\phi(z, \tau) d z,
\end{equation}
where $\Lambda$ is the boundary of the lower half-plane with standard
counter-clockwise orientation, and $\epsilon \to 0^+$. In order to derive 
(\ref{character}), a standard integral representation was used to express the
singular distribution $\delta(x) \theta(x)$, where $\delta, \theta$ are the Dirac
and Heaviside distributions, respectively. 

Equation (\ref{character}) requires a discussion which highlights the physical
aspects of the steady-state selection. At the heart of this discussion lies the 
definition of characteristic function for arbitrary complex values of the parameter
$z$. One way to define $\phi(z), z \in \mathbb{C}$ would be directly through the
integral $\int e^{izx} P(x) dx$. However, we immediately encounter convergence 
issues with this formulation: assume that the {\emph {moment-generating}} function
$m(t) = \int e^{xt} P(x) dx$ has a radius of convergence $t \leq \lambda$. Then 
obviously the characteristic function will also diverge for $\Im (z) < - \lambda$. 
In this way, it is not possible to define a holomorphic function $\phi(z)$ in the whole
complex plane,  unless $\lambda = \infty$. A different generalization of the definition
of $\phi(z)$ is therefore required.

The obvious alternative is to start from the standard characteristic function 
$\phi(z)$ defined on the real axis $\Im (z) = 0$, and to analytically continue it
into the whole complex plane. As we shall see, the convergence radius 
given by the moment-generating function will still play an important role, by determining 
the location and type of singularities of this holomorphic function, but otherwise 
it will be possible to use standard complex analysis techniques to study equation 
(\ref{character}). For example, for a pure exponential distribution of parameter 
$\lambda$, the characteristic function has the form
\begin{equation} \la{steady}
\phi_\infty(z) = \frac{\lambda}{\lambda-iz}.
\end{equation}
Thus, this function has a simple pole at $ z = - i \lambda$, and 
the convergence radius is precisely $\lambda$. However, since $\phi(z)$ 
is now defined everywhere on $\mathbb{C} \setminus \{ - i \lambda \}$ as an analytic function, it is possible to perform the integration in (\ref{character}).  The fact that any pure exponential 
is a steady-state solution, reduces in this language to a simple application of Cauchy's
theorem. 

In this paper, we will consider other distributions related to the exponential, like a
product between a polynomial and an exponential,
\begin{equation}
\label{poly}
p(x,\tau=0)=\lambda\sum_{n=0}^N p_n(0) \frac{(\lambda\, x)^n}{n!}e^{-\lambda\,x},
\end{equation}
or a superposition of exponentials,
\begin{equation} \label{expo}
p(x,\tau)=\int_{\lambda}^\infty d\lambda' \,f(\lambda', \tau)\, \lambda'
e^{-\lambda' x},
\end{equation}
as well as generic sub-asymptotic corrections to the exponential of Gamma-type:
\be  \la{types3}
p(x) \sim \Gamma(\alpha+1, \lambda^{-1}), \quad \alpha > 0. 
\ee
For a distribution of type (\ref{poly}), the analytic continuation of the characteristic function is 
\begin{equation} \la{pole}
\phi(z,\tau) = \sum_{n=0}^N p_n(\tau) \left ( \frac{\lambda}{\lambda - iz}\right )^{n+1},
\end{equation}
while for the type (\ref{expo}), it is
\begin{equation} \la{cut}
\phi(z,\tau) = \int_{\lambda}^{\infty} f(\lambda', \tau)  \frac{\lambda'}{\lambda'-iz} d\lambda'.
\end{equation}
In the first case, the singularities are multiple poles at $z = -i\lambda$, while in the second, 
we have a branch cut extending from $-i\lambda$ to $-i\infty$, with jump function $\lambda'f(\lambda')$
Finally, for the general case (\ref{types3}), the characteristic function combines both 
singularities, having a distribution of simple and multiple poles from $-i\lambda$ to $\infty$. 
It is given by
\begin{equation} \nonumber
\phi(z) = \frac{\Gamma(1-\epsilon)}{\Gamma(1+\alpha)} \left ( -\frac{d}{dz}\right )^n
\frac{1}{\pi}\int_{\lambda}^\infty f(\lambda', \tau)  \frac{\lambda'}{\lambda'-iz} d\lambda',
\end{equation}
where $f(\lambda')$ represents the jump of the function $(z-i\lambda)^{\epsilon-1}$
accross the branch cut, and $\alpha = n- \epsilon,
 0 \le \epsilon < 1, n \in \mathbb{N}, n > 0$.

\section{Direct analysis for some particular types of distributions}

In this section, we illustrate the convergence of coarse-graining for
difference processes towards a pure exponential distribution on a
couple of examples which cover all the types of singularities identified
previously.

\subsection{Polynomial times an exponential}

The integral in (\ref{pp-eq}) has an interesting invariance property. Starting from
an exponential times a polynomial (\ref{poly})
\begin{equation}
p(x,\tau=0)=\lambda\sum_{n=0}^N p_n(0) \frac{(\lambda\, x)^n}{n!}e^{-\lambda\,x},
\end{equation}
with $\sum_{n=0}^N p_n(0)=1$ to ensure proper normalization, the
probability distribution retains the same form
\begin{equation}
p(x,\tau)=\lambda\sum_{n=0}^N p_n(\tau) \frac{(\lambda\,x)^n}{n!}e^{-\lambda\,x}.
\end{equation}
The coefficients $p_n(\tau)$ satisfy a nonlinear evolution equation. We show this 
for the lowest values of $N$.

\begin{itemize}

\item[1)] $N=1$  When the polynomial is linear, the system evolves according to
\begin{equation}
\frac{dp_1}{d\tau}=-\frac{1}{2}p_1^2.
\end{equation}
Therefore, the constant $p_1$ asymptotically decay according to
\begin{equation} \la{rate1}
p_1\simeq 2\tau^{-1}.
\end{equation} 
Thus, the system flows toward the fixed point (\ref{steady-sol}). 

\item[2)] $N=2$ When the polynomial is quadratic, the two independent coefficients
obey
\begin{eqnarray*}
\frac{dp_1}{d\tau}&=&-\frac{1}{2}p_1^2+\frac{1}{2}p_2-\frac{3}{4}p_1p_2-\frac{1}{8}p_2^2\\
\frac{dp_2}{d\tau}&=&-\frac{1}{2}p_1p_2-\frac{3}{4}p_2^2.
\end{eqnarray*}
The last two terms in the first equation and the last term in the
second equation are irrelevant asymptotically and consequently
\begin{equation}
p_1\simeq 4\tau^{-1}\quad p_2\simeq 8\tau^{-2}.
\end{equation}

\item[3)] $N \ge 3$ Generally, the coefficient satisfy
\begin{equation}
\frac{dp_n}{d\tau}=-\frac{1}{2}p_1p_n+\frac{1}{2}p_{n+1}+\cdots
\end{equation}
for $1\leq n\leq N$ with the boundary condition $P_{N+1}=0$.  Here, we
kept only the two asymptotically relevant terms.  The coefficients
decay as follows, $p_n\sim \tau^{-n}$, and it is even possible to obtain
the prefactor
\begin{equation}
p_n(\tau)\simeq \frac{N!}{(N-n)!}2^n\,\tau^{-n}.
\end{equation}

\end{itemize}

Therefore, the asymptotic behavior is independent of the initial
conditions.  Thus, starting from an arbitrary exponential times a
polynomial of the form (\ref{poly}), the system approaches the
purely exponential fixed point
\begin{equation}
p(x)\to \lambda e^{-\lambda\,x}
\end{equation}
as $t\to\infty$.

The average of the variable $x$, $\mathbb{E}(x) =\int_0^\infty dx \, x p(x)$, converges to a constant value according to
\begin{equation}
\label{av-decay}
\mathbb{E}(x)-\lambda^{-1}\simeq 2N\lambda^{-1}\tau^{-1}\simeq 
2N\lambda^{-1}(\ln t)^{-1}.
\end{equation}
Thus, there is a very slow  approach to the steady-state.

\subsection{Laguerre polynomials expansion}

Consider an expansion of the form 

\begin{equation} \label{laguerre}
P_k(x) = \sum_{n=0}^\infty (-1)^n A_n^{(k)} L_n(2x) e^{-x},
\end{equation}
subject to the normalization condition (\ref{normalization}). This expansion is a 
natural reformulation of the previous case, justified by the fact that it corresponds 
to another class of solution for (\ref{steady0}), including non-positive function 
\cite{Bender_Ben-Naim}.
It follows that  the coefficients $A^{(k)}_n$ satisfy the constraint
\begin{equation} \label{condition}
\sum_{n=0}^\infty A^{(k)}_n = 1,
\end{equation}
and that the nonlinear evolution (\ref{pp-eq}) becomes
\begin{equation} \label{system}
A^{(k+1)}_n = \sum_{l=0}^\infty A^{(k)}_l (A^{(k)}_{n+l} + A^{(k)}_{n+l+1}).
\end{equation}
Clearly, if the initial distribution has positive coefficients $\{ A^{(k)}_n \}$, 
they will remain so during the evolution. An asymptotic solution can be found from
(\ref{condition}) and (\ref{system}), asuming that 
\begin{equation}
A_0 \to 1, \quad A_{n \ge 1} \to 0, \quad A_1 \ge A_{k \ge 2}.
\end{equation} 
Under these assumptions, the long-time evolution equations become:
\begin{equation}
\begin{array}{lcl}
\frac{d A_1}{d t} & = & A_2 - A_1^2 + O(A_1A_2), \\
&& \\
\frac{d A_k}{d t} & = & A_{k+1} - A_k A_1 + O(A_k A_2), \,\, k \ge 2,
\end{array}
\end{equation}
which has the solution
\begin{equation}
A_1(\tau) = \frac{\alpha}{\tau} + O(\tau^{-2}), 
\end{equation}
\begin{equation}
A_k(\tau) = \frac{\Gamma(\alpha+1)}{\Gamma(\alpha-k+1)}\tau^{-k} + O(\tau^{-k-1}).
\end{equation}
Clearly, if we require that $A_k \ge 0$ throughout the evolution, it follows that
starting some index $k=N$, we must have $A_{k \ge N} = 0$, and that 
\begin{equation}
A_1 = \frac{N}{\tau}, \quad A_k = \frac{N!}{(N-k)!}\tau^{-k}, \,\, 1< k < N.
\end{equation}
The positivity condition $A_k \ge 0$ has a clear physical significance: if it 
is not satisfied, we see that starting some value of $k$, all coefficients become
negative. Therefore, at fixed $\tau$, the tail of the distribution $P_k(x)$ is not in the exponential class, in fact it can be very far from exponential. The exponential-type asymptote is guaranteed only by the positivity condition, which effectively truncates the expansion (\ref{laguerre}). In turn, this condition alone is enough to determine the asymptotic behavior of the coefficients $A_k$ and therefore of the average $\mathbb{E}[x] = 1+ \frac{2N}{\tau}$ (compare with (\ref{av-decay})). 

\subsection{Sum of exponentials}

The integral in (\ref{pp-eq}) has another important property. Let us
start with a sum of $N+1$ exponentials
\begin{equation}
\label{sum-exp-0} p(x,\tau=0)=\sum_{j=0}^N A_j(0)\lambda_j
e^{-\lambda_j x},
\end{equation}
with monotonic decay constants $\lambda_j<\lambda_{j+1}$. Then, the
solution remains a sum of exponentials
\begin{equation}
\label{sum-exp} p(x,\tau)=\sum_{j=0}^N A_j(\tau)\lambda_j
e^{-\lambda_j x}.
\end{equation}
The coefficients satisfy $\sum_{j=0}^N A_j=1$ to guarantee
normalization. Substituting (\ref{sum-exp}) into the master equation
(\ref{p-eq}), the coefficients evolve according to
\begin{equation}
\frac{dA_k}{d\tau}=-A_k+2A_k\sum_{j=0}^N
A_j\frac{\lambda_j}{\lambda_j+\lambda_k}.
\end{equation}
Replacing the term $-A_k$ with $-A_k\sum_{j=0}^N A_j$, we rewrite
this equation.
\begin{equation}
\frac{dA_k}{d\tau}=A_k\sum_{j=0}^N
A_j\frac{\lambda_j-\lambda_k}{\lambda_j+\lambda_k}.
\end{equation}
There are $N+1$ steady-states of the type $A_k=\delta_{k,j}$ but only
the steady-state $A_k=\delta_{k,0}$ is stable. Small perturbations
to this steady-state decay exponentially fast
\begin{equation}
\label{ak-decay} A_k\sim
\exp\left[-\frac{\lambda_k-\lambda_0}{\lambda_k+\lambda_0}\,\tau\right].
\end{equation}
It is simple to show recursively (first for $k=N$ and then for
$k=N-1$ and so on) that the coefficients $A_k$ decay exponentially
as in (\ref{ak-decay}). Thus, starting from a sum of exponentials
(\ref{sum-exp-0}), the steady state (\ref{steady-sol}) is selected,
\begin{equation}
\label{select} p(x)\to \lambda_0\exp(-\lambda_0\,x)
\end{equation}
as $\tau\to\infty$.

The approach toward the steady state is dominated by the coefficient
$A_1$:
\begin{equation}
\mathbb{E}[x] -\lambda_0^{-1}\sim
\exp\left[-\frac{\lambda_1-\lambda_0}{\lambda_1+\lambda_0}\,\tau\right]
\sim t^{-\frac{\lambda_1-\lambda_0}{\lambda_1+\lambda_0}}.
\end{equation}

\subsection{Integral of exponentials}

It is straightforward to generalize the above from a sum to an
integral
\begin{equation} 
p(x,\tau)=\int_{\lambda_0}^\infty d\lambda \,f(\lambda, \tau)\, \lambda
e^{-\lambda x}.
\end{equation}
with $\int_{\lambda_0}^\infty\, d\lambda\, f(\lambda)=1$. Then, the
coefficient $f(\lambda)$ evolves according to
\begin{equation}
\frac{\partial f(\lambda)}{\partial
\tau}=f(\lambda)\,\int_{\lambda_0}^\infty d\lambda'\,
f(\lambda')\,\frac{\lambda'-\lambda}{\lambda'+\lambda}.
\end{equation}
The only stable steady state is $f(\lambda)=\delta(\lambda_0)$ as in
(\ref{select}) and for all $\lambda>\lambda_0$, the coefficients
$f(\lambda)$ decays exponentially,
\begin{equation}
f(\lambda)\sim
\exp\left[-\frac{\lambda-\lambda_0}{\lambda+\lambda_0}\,\tau\right].
\end{equation}

\section{General analysis of convergence to steady-state}

The conclusions of previous sections show that starting from a characteristic function
featuring a multiple pole (\ref{pole}) at $z = -i\lambda$, (\ref{pp-eq}) leads to the 
steady-state (\ref{steady}), with power-law convergence of the first moment. The case
of a branch cut extending from $-i\lambda_0$ to $-i\infty$ (\ref{cut}) leads to 
(\ref{steady}) with exponentially fast convergence. These two complementary situations 
allow to explain the general steady-state selection process.

The generic situation is covered by the the class of distributions (\ref{types3}), where 
$\alpha = n- \epsilon, 0 \le \epsilon < 1, n \in \mathbb{N}, n > 0$. The characteristic function is 
\begin{equation} \la{chargam}
\phi(z) = \left ( \frac{\lambda}{\lambda-iz}\right )^{\alpha+1},
\end{equation}
so it has the mixted representation
\begin{equation} \nonumber
\phi(z) = \frac{\Gamma(1-\epsilon)}{\Gamma(1+\alpha)} \left ( -\frac{d}{dz}\right )^n
\frac{1}{\pi}\int_{\lambda}^\infty f(\lambda', \tau)  \frac{\lambda'}{\lambda'-iz} d\lambda',
\end{equation}
where $f(\lambda')$ represents the jump of the function $(z-i\lambda)^{\epsilon-1}$
accross the branch cut. This can be interpreted as a combination of distributions of 
type (\ref{pole}), where now the parameter of the exponential ranges from $\lambda$ to 
$\infty$. By the same argument as before, the combination of exponentials decays exponentially fast towards the lowest value of $\lambda$, while the pole of order $n+1$ decays algebraically
to the pure exponential. This separation of time scales in the convergence process allows to
conclude that the steady state corresponding to an initial choice of type (\ref{types3}) 
is $\lambda e^{-\lambda x}$, and the first moment converges to $\lambda^{-1}$ as 
\begin{equation}
\mathbb{E} [x] -\lambda^{-1} = \frac{2n}{\lambda \tau}.
\end{equation} 

\section{Conclusions}
Coarse graining of difference processes represents a fundamental generalization of 
standard coarse graining with averaging. This procedure is a natural description for 
relevant phenomena ranging from multi-species stochastic processes to socio-economics.
In this paper, we have identified the class of stead-states for this process, and shown how
a particular, universal limiting distribution is chosen, as well as the convergence rate towards the
steady-state.  

\section*{Acknowledgements}

The author is grateful to C. Bender for introducing him to the problem and for useful discussions.

\end{document}